\documentclass[]{interact}

\usepackage{epstopdf}
\usepackage[caption=false]{subfig}
\usepackage{apacite}
\usepackage[longnamesfirst,sort]{natbib}
\bibpunct[, ]{(}{)}{;}{a}{,}{,}
\renewcommand\bibfont{\fontsize{10}{12}\selectfont}

\theoremstyle{plain}
\newtheorem{theorem}{Theorem}[section]
\newtheorem{lemma}[theorem]{Lemma}

{
	\newtheorem{assumption}{Assumption}
	
}
\theoremstyle{definition}

\usepackage{color, xcolor}
\usepackage{hypernat}
\usepackage[colorlinks]{hyperref}
\hypersetup{citecolor=blue}
\theoremstyle{remark}
\newtheorem{remark}{Remark}

\begin{document}

\articletype{ARTICLE TEMPLATE}

\title{Distributed optimization problem with communication delays for heterogeneous linear multi-agent systems}

\author{
\name{Farshad Rahimi\textsuperscript{a}\thanks{CONTACT Farshad Rahimi. Email: fa\_rahimi@sut.ac.ir}\textsuperscript{a}}
\affil{\textsuperscript{a}Department of Electrical Engineering, Sahand University of Technology,  Tabriz, Iran.
	Tel.: +98-83-46422284.}
}

\maketitle

\begin{abstract}
This work studies the problem of distributed optimization in heterogeneous linear multi-agent systems. Instead of relying on a perfect communication network as in many existing distributed optimization approaches, we considered two important issues related to communication networks. First, assumed that communication delays exist when each agent receives information from its neighbors. Further, since communications networks are generally unreliable, we assumed each agent interacted with the other agents through random digraphs. Finally, to prove convergence to optimal solutions, delay-dependent sufficient conditions are derived in the form of linear matrix inequality. An analysis of numerical simulation results is presented to demonstrate the effectiveness of the introduced approach.
\end{abstract}

\section{Introduction}
\label{intro}
In distributed optimization problem for multi-agent systems, the sum of multiple functions is minimized or maximized by individuals working together to achieve the objective function. This type of problem usually involves each agent having access to only one local objective function and exchanging information with its neighbors \cite{f3}. In science and engineering, a variety of applications arise from distributed optimization using multi-agent systems, such as resource allocation using networks \cite{f4}, energy and thermal comfort optimization in smart building \cite{f5}, formation control \cite{f6,f24,f28}, etc. For instance, in reference \cite{f23}, an economic dispatch method with constrained distribution is suggested to deal with transport line losses and total load demand. There has been a great deal of research into distributed optimization using multi-agent systems. Many of these works are devoted to discrete-time formulations \cite{f9,f10,f29,fn55}. The readers interested in the distributed optimization problem can refer to \cite{f27}, in which models and algorithms with a variety of constraints have been examined in relation to the model structure, the algorithm type, and the communication topology.
In recent years, the study of continuous-time distributed optimization has received much attention in part due to its simple application for convergence analysis and also for its potential applications in continuous-time physical systems \cite{cm1,f7,f8}. For instance, researchers in \cite{f11} investigated the problem of distributed optimization in continuous-time multiagent systems with parametric uncertainties.\\
The existing literature on the distributed optimization problem mostly deals with homogeneous multi-agent systems, that is, all systems have the same dynamics. Therefore, first motivation in this paper is about consider the problem of distributed optimization for heterogenous multi-agent systems. The communication network is required for each agent within a multi-agent system in order to receive and transmit information.  The communication medium usually has a finite bandwidth, and the agents have a limited capacity for computing and sending their results \cite{f12}, thus, during data exchange between agents in a network environment, there are inevitable delays. Multi-agent systems suffer from performance and stability degradation when there is a communication delay.  In the establishment of stability conditions for a networked system, it is crucial to consider time delays \cite{f16,f13,f14,f15}. With the help of scattering transformation, a class of distributed optimization problems subject to inter-agent communication delays has been introduced \cite{f13}, where a convergence property of the problem can be guaranteed with the integration of inter-agent communication delays. A time-varying communication delay is also presented in \cite{f14} for the distributed optimization problem of the continuous-time multi-agent systems. By designing a continuous-time algorithm with small and large delays, the authors in \cite{f15} investigated the distributed constrained optimization problem. To solve this problem, they used switching techniques to develop a switched delay algorithm.\\
Concerning the exchange of data among agents, the switching topology problem is an important  issue in controller design for multi-agent systems.  Because communication networks are vulnerable, the structure of communication graphs in multiagent systems can be changed by a variety of causes \cite{f22}. There are three types of topological graphs in communication networks: random graphs, directed connected graphs, and undirected connected graphs. The random graphs approach lends more credit to the proposed approaches since real communication topology among agents can be changed during the implementation of the controllers. A random graph can display this change, but a directed graph or an undirected graph cannot. In recent years, some works on random graphs in the distributed optimization problem have been published in \cite{f20,f21}. The authors in the work \cite{f20} examined a collaborative optimization problem involving a sum of convex functions under the assumption that agents use local information and random topologies in their interconnections. To solve a distributed optimization problem over Markovian switching random network communication, the work \cite{f21} proposed a consensus subgradient algorithm where unlike most existing algorithms, that approach had two different time-scale step-sizes. The authors of \cite{f25} by considering the fact that sometimes random digraphs would be caused by unreliable networks  and Denial-of-Service attacks, investigated the problem of distributed convert optimization for linear multi-agent systems.
\\
\textbf{The main contributions and the closed works:} Motivated by the aforementioned papers and discussions, this paper is devoted to introducing a distributed optimization approach for  continuous-time heterogenous multi-agent system with time-varying communication delays. The formulation of our optimization problem is based on \cite{f30,f25} on which they introduced a distributed optimization problem for  continuous-time multi-agent system. As opposed to references \cite{f30,f25}, this paper considers communication delays in formulating the distributed optimization problem for networked heterogeneous systems, then establishes a stability condition for estimating the upper bound on the time delay. Finally, we derived  the relationship between optimal solutions and the estimation of  the upper bound of time delay  for  the distributed optimization of heterogenous multi-agent systems. The proposed approach uses a gradient-based consensus strategy with a switching system method to tackle the distributed optimization problem of a heterogeneous linear multi-agent system.  The aforementioned method can seek optimal solutions while dealing with communication delays over random digraphs. The main contributions of this study can be highlighted as follows: \\
(1)	A distributed optimization approach is established for heterogenous linear multi-agent systems where each agent receives the its requirement information from its neighbors with communication delays over random digraphs.\\
(2)	A Lyapunov-Krasovskii stability argument is used to verify the convergence of the proposed distributed optimization scheme. In addition, a delay-dependent sufficient condition based on a LMI is derived to find the maximum delay which the considered optimization approach can theoretically tolerate.  \\
Following is an overview of how this paper is organized. Graph theory and convex analysis are introduced in Section \ref{sec1} with preliminaries on basic notation and related concepts. After that, an optimization problem for multi-agent systems with heterogeneous linear behavior is formulated. The results of the study are outlined in Section \ref{sec_mr}.  Section \ref{sec_sr} gives numerical examples of the proposed algorithm that illustrate its effectiveness. Finally, Section \ref{sec_c} concludes this article.\\
\textbf{Notations:} The notations used in this paper are fairly standard.
Denote $\mathbb{R}$, ${{\mathbb{R}}^{p}}$ and ${{\mathbb{R}}^{p\times q}}$ as the sets of the real numbers, real $p$-dimensional vectors and real $p\times q$  matrices, respectively. We write $B>0$ and $B\ge 0$ to indicate that $B$ is a positive definite and a positive semidefinite, respectively. $diag\{{{b}_{1}},...,{{b}_{n}}\}$ means a diagonal matrix with the entries ${{b}_{j}},\,\,j=1,2,...,n$.

\section{Preliminaries and Formulation Problem}
\label{sec1}
\subsection{Graphs}
Here,  we give some preliminaries concerned with graph theory. The notation ${\Xi }\left( \mathcal{V},\mathcal{E},\Im  \right)$ demonstrates a graph where $\mathcal{V}=\left\{ 1,\ldots ,N \right\}$ is a set of nodes. A set of edges can be defined as $\mathcal{E}\subseteq \mathcal{V}\times \mathcal{V}$. $N$ describes the numbers of agents that exists in the network.  The neighborhood set of vertex $i$ is described by ${{{N}}_{{i}}}=\left\{ j\text{ }\!\!|\!\!\text{ }\left( i,j \right)\in \mathcal{E} \right\}$. An unilateral path from node $i$ to $j$  means agent $i$ can send its information to agent $j$. $A=\left[ {{a}_{ij}} \right]\in {{R}^{N\times N}}$ is an adjacency matrix of the graph ${\Xi }\left( \mathcal{V},\mathcal{E} \right)$.  The information flow in ${\Xi }$ is described by an entire ${a_{ij}}$ satisfying ${{a}_{ij}}=1$ if   $i \ne j $ and $(i,j)\in \mathcal{E}$, otherwise ${{a}_{ij}}=0$. The Laplace matrix of graph ${\Xi }$ is given by $L=H-A$, where $H=diag\{{{h}_{i}}\}\in {{R}^{N\times N}}$ be  an in-degree matrix with ${{h}_{i}}=\sum\limits_{j\in {{N}_{i}}}{{{a}_{ij}}}$.
Since, Markovian random digraph would be employed in the structure of the proposed approach, it is necessary to definite them at the beginning.\\
Let $\Omega (t)=\{\mathcal{V},{{\mathcal{E}}_{r(t)}}\}$ be a time varying digraph with ${{\mathcal{E}}_{r(t)}}$ being a set of edges, and $r(t):[0,\infty )\to S=\{1,2,...,s\}$ is a piecewise constant function with   being an index set of possible digraphs. The piecewise-constant function $r(t)$ is a Markovian signal. $A(t)=\left[ a_{ij}^{r(t)} \right]$ is the adjacency  matrix, where $a_{ij}^{r(t)}>0$ if $(i,j)\in {{\mathcal{E}}_{r(t)}},$ else   $a_{ij}^{r(t)}=0$. The neighboring set is denoted by ${{N}_{i}}(t)=\left\{ j\in \mathcal{V},\,(j,i)\in {{\mathcal{E}}_{r(t)}} \right\}$ . Denote $L(t)=\left[ l_{ij}^{r(t)} \right]$ , where $l_{ii}^{r(t)}=\sum\limits_{j=1}^{N}{a_{ij}^{r(t)}}$ and $l_{ij}^{r(t)}=-a_{ij}^{r(t)},i\ne j$.

\subsection{Problem Formulation}
Consider a multi-agent system with $N$ agents and a graph ${\Xi }$ describing their interactions. Agent $i$ operates according to the following dynamic:
\begin{align}
	& {{{\dot{x}}}_{i}}(t)={{A}_{i}}{{x}_{i}}(t)+{{B}_{i}}{{u}_{i}}(t), \cr
	& {{y}_{i}}(t)={{C}_{i}}{{x}_{i}}(t),\,\,\,i\in \mathcal{V}
	\label{eq1}
\end{align}

where ${{x}_{i}}(t)\in {{\mathbb{R}}^{{{n}_{i}}}}$ stands for the state of agent $i$, ${{u}_{i}}(t)\in {{\mathbb{R}}^{{{p}_{i}}}}$ represents the control protocol of agent $i$, the output is described by ${{y}_{i}}(t)\in {{\mathbb{R}}^{{{q}_{i}}}}$. Matrices ${{A}_{i}}\in {{\mathbb{R}}^{{{n}_{i}}\times {{n}_{i}}}}$, ${{B}_{i}}\in {{\mathbb{R}}^{{{n}_{i}}\times {{p}_{i}}}}$ , and ${{C}_{i}}\in {{\mathbb{R}}^{{{q}_{i}}\times {{n}_{i}}}}$ are known matrices with appropriate dimensions.  As a prerequisite for determining the stability of the proposed approach, we need the following assumption.
\begin{assumption}
 It is assumed that system (\ref{eq1}) is controllable, meaning:
$rank\left[ \begin{matrix}
	{{C}_{i}}{{B}_{i}} & {{0}_{{{q}_{i}}\times {{p}_{i}}}}  \\
	-{{A}_{i}}{{B}_{i}} & {{B}_{i}}  \\
\end{matrix} \right]={{n}_{i}}+{{q}_{i}},\,\,\,\,i\in \mathcal{V}.$
\label{asp1}
	\end{assumption}
\subsection{Control Objective}
This paper intends to design a distributed control protocol for the system (\ref{eq1})  such that the output all agents converges to the optimal state.  Through a communication network, each agent utilizes the information of its neighbors through its control input ${{u}_{i}}(t)$. Agent $i$ receives the information from its neighbors with a time-varying delay. We can formulate our distributed optimization problem as follows:
\begin{equation}
F(\theta )=\sum\limits_{i=1}^{N}{{{f}_{i}}(\theta )},\,\,\,\,\,\theta \in {{\mathbb{R}}^{{{q}_{i}}}}.
	\label{eq2}
\end{equation}
where ${{f}_{i}}(\theta ):\,{{\mathbb{R}}^{{{q}_{i}}}}\to \mathbb{R}$ represents the local cost function that is only known to agent $i$. $\theta $ is the global decision variable that must be optimized. As an equivalent, we can solve problem (\ref{eq2}) as follows:
\begin{equation}
\begin{matrix}
	\begin{matrix}
		\underset{y\in {{\mathbb{R}}^{{N{q}_{i}}}}}{\mathop{\min }}\, & \tilde{f}(y)=  \\
	\end{matrix}\sum\limits_{i=1}^{N}{{{f}_{i}}({y}_{i} )},\,\,\,\,\,{{y}_{i}}\in {{\mathbb{R}}^{{{q}_{i}}}}  \\
	s.t.\,\,(1)\,\,and\,\,{{y}_{i}}={{y}_{j}},\,\,\forall i,j\in \mathcal{V}=\left\{ 1,\ldots ,N \right\}.  \\
\end{matrix}
	\label{eq3}
\end{equation}
The local estimate ${{y}_{i}}\in {{\mathbb{R}}^{{{q}_{i}}}}$ here is the optimal solution$\,{{\theta }^{*}}$. Additionally, all estimates ${{y}_{i}},\,\,i\in \left\{ 1,...,N \right\}$ can be augmented by $y=col({{y}_{1}},....,{{y}_{N}})$.
\\
There are some assumptions that need to be made before solving the optimization problem (\ref{eq3}).
\begin{assumption}
 It is possible to find ${{y}^{*}}={{1}_{{{q}_{i}}}}\otimes {{\theta }^{*}}$ which minimizes the team cost function, i.e., $\tilde{f}({{y}^{*}})={{\min }_{\theta \in {{\mathbb{R}}^{{{q}_{i}}}}}}F(\theta )$.
 \label{ass1}
\end{assumption}
\begin{assumption}
 Suppose that all functions ${{f}_{i}}:{{\mathbb{R}}^{{{q}_{i}}}}\to \mathbb{R},\,\,i\in \left\{ 1,...,N \right\}$ are differentiable, strongly convex, and their gradients locally Lipschitz on ${{\mathbb{R}}^{{{q}_{i}}}}$.
 \label{ass2}
\end{assumption}
Locally Lipschitz means $\left\| \nabla {{f}_{i}}({{x}_{i}})-\nabla {{f}_{i}}({{y}_{i}}) \right\|\le {{l}_{i}}\left\| \nabla {{x}_{i}}-\nabla {{y}_{i}} \right\|,\,\,\forall {{x}_{i}},{{y}_{i}}\in {{\mathbb{R}}^{{{q}_{i}}}}$, where $\nabla {{f}_{i}}({{x}_{i}})$ and $\nabla {{f}_{i}}({{y}_{i}})$ are the gradients, and ${{l}_{i}}>0$ is the Lipschitz constant.\\
\label{att1}
Due to physical uncertainties, such as failure and packet loss during digital communication, the communication networks between multiple agents are not reliable. Furthermore, multi-agent systems may also be viewed as unreliable because they are vulnerable to attack by the adversary due to the vulnerability of the communication networks during the exchange of information between agents. This paper considers an unreliable network consisting of $N$ agents where communication links are time-varying and fail-prone with certain probabilities.  This feature is taken into consideration in \cite{f31} by using a random Markov chain model.  We will assume that $r(t)$ is a right-continuous Markov-vian process on the probability space taking values from a finite state space $S=\left\{1,2,...,s\right\}$ with an infinitesimal generator $\Upsilon =({\gamma}_{pq})$, given by ${{P}_{pq}}(t)=\Pr ob\left\{ r(t+h)=q|r(t)=p \right\}={{\gamma }_{pq}}h+o(h),$ if $p\ne q$ , else, $1+{{\gamma }_{pp}}h+o(h),$ where ${{\gamma }_{pq}}\ge 0$ is the transition rate from the state $p$ to the state $q$ , while ${{\gamma }_{pp}}=-\sum\limits_{q=1,p\ne q}^{s}{{{\gamma }_{pq}}}$ , and $o(h)$ satisfies: $\underset{h\to 0}{\mathop{\lim }}\,o(t)/h=0$ .
\subsection{Optimization Formulation}
A key objective is to implement a distributed optimization algorithm, ${u}_{i}(t)$, so that the output of each agent cooperatively seeks to optimal ${\theta}^*$ over random digraphs. This is the reformulated problem of (\ref{eq3})
\begin{align}
	\begin{matrix}
		\begin{matrix}
			\underset{y\in {{\mathbb{R}}^{{{q}_{i}}}}}{\mathop{\min }}\, & \tilde{f}(y)=  \\
		\end{matrix}\sum\limits_{i=1}^{N}{{{f}_{i}}({{y}_{i}})},\,\,\,\,\,{{y}_{i}}\in {{\mathbb{R}}^{{{q}_{i}}}}  \\
		s.t.\,\,(1)\,\,and\,\,(L(t)\otimes {{I}_{q}})y=0.  \\
	\end{matrix}
\label{eq4}
\end{align}
Due to unreliable networks described in the subsection \ref{att1}, the underlying topologies are time-varying and random. As stated in \cite{f31}, define that a sequence of Laplacian matrices $\left\{ L(t) \right\}$ admits a common stationary distribution $\pi >0$ if $L(t)\pi =0$. Let ${{L}_{s}}(t)$ be a mirror of ${{L}_{un}}(t),$ i.e.,  ${{L}_{s}}(t)=({{L}_{un}}(t)+L_{un}^{T}(t))/2,$where ${{L}_{un}}=\sum\limits_{p=1}^{N}{{{L}_{p}}(t)}$ is Laplacian matrix of a union of digraphs. Define the minimum cut of $\left\{ L(t) \right\}$ as ${{l}_{c}}(t)={{\min }_{S\subset V,S\ne 0}}\sum\limits_{i\in S,j\in \bar{S}}{{{L}_{s}}(t)}$, where $\bar{S}$ is the complement  of $S$. Then, we say that this sequence of $\left\{ L(t) \right\}$ has a minimum cut $c$, if ${{l}_{c}}(t)\ge c>0$ .
\begin{assumption}
If the set of Laplacian matrices  $\left\{ L(t) \right\}$ has a stationary distribution $\pi $ with a minimum cut, then we think of the sequence as having a stationary distribution.
\label{ass6}
\end{assumption}
\begin{lemma}
	\label{lem1}
	 By Assumption \ref{ass6}, let  $\left\{ L(t) \right\}$ be the Laplacian matrix with a stationary distribution $\pi>0 $. Then, there exists a weighted matrix $Q(t)=L(t)\Pi +\Pi {{L}^{T}}(t),\,\,\Pi =diag\left\{ \pi  \right\}$ , so that for ${{\pi }_{\min }}={{\min }_{p\in V}}\left\{ {{\pi }_{p}} \right\}$ and a vector $\zeta \in {{\mathbb{R}}^{N}}$ satisfying ${{\pi }^{T}}\zeta =0$ ,  we have
	 \begin{align}
	 	{{\zeta }^{T}}Q(t)\zeta =\sum\limits_{i=1}^{N}{\sum\limits_{j\in N(t)}{{{Q}_{ij}}(t){{({{\zeta }_{i}}-{{\zeta }_{j}})}^{2}}\ge \frac{{{\pi }_{\min }}c}{{{N}^{2}}}}{{\left\| \zeta  \right\|}^{2}}}
	 	\label{eq5}
	 \end{align}
\end{lemma}

\subsection{Communication Delays}
The paper also considers an important issue in unreliable networks which are common in multi-agent systems, i.e.,  communication delays. It is not immediately possible to access neighbor information in agent $i$, i.e., at step time $t$, the amount of ${{y}_{j}}(t-{{d}_{i}}(t))$ is available instead of ${{y}_{j}}(t)$. We assume that each agent in the networked system receives information from its neighbors with a time-varying delay ${{d}_{i}}(t)$, which it is bounded as $0\le {{d}_{i}}(t)\le \bar{{{d}_{i}}},\,\,\,\dot{{{d}_{i}}}(t)\le {{\varpi }_{i}} $. In whats follows, we need to implement the effects of communication delay in the process of establishing the stability. To do this, we will use the following lemma.
\begin{lemma}
	For any symmetric and positive definite matrix ${{R}_{1}}$ scalar $0\le {{d }_{i}}(t)\le {{{d }}_{m}}$ and  vector function $\dot{Z}(t):[-{{d }_{m}},0]$ making the following inequality is defined, it holds that
	\begin{align*}
		&{{d }_{m}}\int\limits_{t-{{d }_{m}}}^{t}{{{{\dot{Z}}}^{T}}\left( s \right)}{{R}_{1}}\dot{Z}\left( s \right)ds\le {{\left[ \begin{matrix}
					Z(t)  \\
					Z(t-{{d }_{m}})  \\
					Z(t-d (t))  \\
				\end{matrix} \right]}^{T}} \cr
		& \,\,\,\,\,\,\,\,\,\times \left[ \begin{matrix}
			-{{R}_{1}} & 0 & {{R}_{1}}  \\
			* & -{{R}_{1}} & {{R}_{1}}  \\
			* & * & -2{{R}_{1}}  \\
		\end{matrix} \right]\left[ \begin{matrix}
			Z(t)  \\
			Z(t-{{d }_{m}})  \\
			Z(t-d (t))  \\
		\end{matrix} \right] 	.
	\end{align*}
	\label{lem11}	
\end{lemma}
In order to analyze the effects of communication delay over random digraphs, we propose the following distributed algorithm.
\begin{align}
	{{u}_{i}}(t)=-{{K}_{i}}{{x}_{i}}(t)-({{U}_{i}}-{{K}_{i}}{{X}_{i}}){{\eta }_{i}}(t)+{{W}_{i}}{{v}_{i}}(t),
	\label{eq6}
\end{align}
\begin{align}
    & {{v}_{i}}(t)=\nabla {{f}_{i}}({{y}_{i}})-\beta {{e}_{\eta {{z}_{i}}}}-\alpha \beta {{e}_{{{y}_{i}}}}, \label{eq7}\\
	& {{{\dot{\eta }}}_{i}}(t)={{v}_{i}}(t),\label{eq8}\\
	& {{{\dot{z}}}_{i}}(t)=\alpha \beta {{e}_{{{y}_{i}}}},\,\, \label{eq9} \,\,\,i\in V
\end{align}
in which there are positive constant gains $\alpha \beta \in \mathbb{R}$, ${{v}_{i}}$ stands for the intermediate variable. The parameters ${{e}_{\eta {{z}_{i}}}},\,\,{{e}_{{{y}_{i}}}}$ represents consensus errors
\begin{align}
	\resizebox{.95\hsize}{!}{$
{{e}_{\eta {{z}_{i}}}}=	\sum\limits_{j\in {{N}_{i}}(t)}{{{a}_{ij}}(t)\Big({{\eta }_{i}}(t-{{d}_{i} }(t))-{{\eta }_{j}}(t-{{d}_{i} }(t))+{{z}_{i}}(t-{{d}_{i}}(t))-{{z}_{j}}(t-{{d}_{i} }(t))\Big)}$},
\label{eq10}
\end{align}
\begin{align}
{{e}_{{{y}_{i}}}}=\sum\limits_{j\in {{N}_{i}}(t)}{{{a}_{ij}}(t)\Big({{y}_{i}}(t-{{d}_{i} }(t))-{{y}_{j}}(t-{{d}_{i}}(t))\Big)},
\label{eq11}
\end{align}
in which auxiliary states are supposed to be denoted by ${{\eta }_{i}}(t),{{z}_{i}}(t)\in {{\mathbb{R}}^{q}}$ . The following equations have a solution for ${{K}_{i}}\in {{\mathbb{R}}^{{{p}_{i}}\times {{n}_{i}}}}$ and $({{U}_{i}},{{W}_{i}},{{X}_{i}})$.
\begin{align}
	{{B}_{i}}{{U}_{i}}={{A}_{i}}{{X}_{i}},\,\,{{B}_{i}}{{W}_{i}}={{X}_{i}},\,\,{{C}_{i}}{{X}_{i}}={{I}_{q}}\,,\,\,\,i\in V.
	\label{eq12}
\end{align}
Based on our assumptions, we assumed $\tilde{L}(t)=L(t)\otimes {{I}_{q}}.$ When (\ref{eq6})-(\ref{eq11}) are substituted into (\ref{eq1}), the following closed-loop systems are produced.
\begin{align}
		&	\resizebox{.98\hsize}{!}{$ \dot{X}=\left[ \begin{matrix}
			\dot{x}(t)  \\
			\dot{\eta }(t)  \\
			\dot{z}(t)  \\
		\end{matrix} \right]=\left[ \begin{matrix}
			(A-BK)x(t)+BWv(t)-B(U-KX)\eta (t)   \cr
			-\nabla \tilde{f}(y(t))-\alpha \beta \tilde{L}(t)y(t-d(t))-\beta \tilde{L}(t)\Big(\eta (t-d(t)) +z(t-d(t))\Big)  \cr
			\alpha \beta \tilde{L}(t)y(t-d(t))  \cr
		\end{matrix} \right] $} \cr &
		 \dot{\eta }(t)=v(t),
		\label{eq13}
\end{align}
in which
\begin{align*}
		& x(t)=diag\left[ {{x}_{1}}(t),\,{{x}_{2}}(t),...,{{x}_{N}}(t) \right],\,\, \\ & \eta (t)=diag\left[ {{\eta }_{1}}(t),\,{{\eta }_{1}}(t),...,{{\eta }_{N}}(t) \right], \\
		& v(t)=diag\left[ {{v}_{1}}(t),\,{{v}_{2}}(t),...,{{v}_{N}}(t) \right],\,\,  \\ &  z(t)=diag\left[ {{z}_{1}}(t),{{z}_{2}}(t),...,{{z}_{N}}(t) \right],
		\\ &  d(t)=diag\left[ {{d}_{1}}(t),{{d}_{2}}(t),...,{{d}_{N}}(t) \right], \\
		& \nabla \tilde{f}(y(t))=diag\left[ \nabla {{f}_{1}}({{y}_{1}}(t)),\nabla {{f}_{2}}({{y}_{2}}(t)),...,\nabla {{f}_{N}}({{y}_{N}}(t)) \right],\,\, \\
		& A=diag\left[ {{A}_{1}},\,{{A}_{2}},...,{{A}_{N}} \right],\,\,B=diag\left[ {{B}_{1}},{{B}_{2}},...,{{B}_{N}} \right], \\
		& C=diag\left[ {{C}_{1}},\,{{C}_{2}},...,{{C}_{N}} \right],\,\,K=diag\left[ {{K}_{1}},{{K}_{2}},...,{{K}_{N}} \right], \\
		& U=diag\left[ {{U}_{1}},\,{{U}_{2}},...,{{U}_{N}} \right],\,\,W=diag\left[ {{W}_{1}},{{W}_{2}},...,{{W}_{N}} \right], \\
		& X=diag\left[ {{X}_{1}},{{X}_{2}},...,{{X}_{N}} \right].
\end{align*}
\section{Main Results}
\label{sec_mr}
In the sequel, we present our results for the distributed optimization problem for heterogenous linear multi-agent systems in the presence of communication delays over random digraphs. Delay-dependent conditions are presented below. In the following Theorem, sufﬁcient conditions are provided to ascertain the convergence of the proposed distributed optimization scheme to the optimal solution.
\begin{theorem}
\label{thm1}
		We assumed  Assumptions \ref{ass1}-\ref{ass2} are satisfied, under the introduced distributed optimization approach in (\ref{eq6})-(\ref{eq11}),  system (\ref{eq4}) with initial value in ${{x}_{i}}(0),{{\eta }_{i}}(0)$ and ${{z}_{i}}(0)$  convergent to the optimal solution for any delay $d(t)\in \left[ 0,\bar{d} \right]$  with $\dot{d}(t)\le \varpi $ and $\varpi \in \left[ 0,1 \right)$, if ${{K}_{i}}$ is selected to that ${{A}_{i}}-{{B}_{i}}{{K}_{i}}$ is Hurwitz, ${{U}_{i}},\,{{W}_{i}},\,{{X}_{i}}$ are solutions to (\ref{eq12}), and if there exist  positive deﬁnite matrices  ${{P}_{a}},{{P}_{1}},{{P}_{2}},{{P}_{3}},{{P}_{4}},{{Q}_{1}},{{Q}_{2}},{{Q}_{3}}$ such that the following LMI are feasible:
		\begin{align}
			\Pi =\left[ \begin{matrix}
				{{\Pi }_{11}} & 0 & 0 & {{\Pi }_{14}} & 0 & 0 & {{\Pi }_{17}} & {{\Pi }_{18}}  \\
				* & {{\Pi }_{22}} & 0 & {{\Pi }_{24}} & 0 & 0 & 0 & {{\Pi }_{28}} \\
				* & * & {{\Pi }_{33}} & {{\Pi }_{34}} & 0 & 0 & 0 & 0  \\
				* & * & * & {{\Pi }_{44}} & {{\Pi }_{45}} & 0 & {{\Pi }_{47}} & {{\Pi }_{48}}  \\
				* & * & * & * & {{\Pi }_{55}} & 0 & {{\Pi }_{57}} & {{\Pi }_{58}}  \\
				* & * & * & * & * & {{\Pi }_{66}} & {{\Pi }_{67}} & 0  \\
				* & * & * & * & * & * & {{\Pi }_{77}} & {{\Pi }_{78}}  \\
				* & * & * & * & * & * & * & {{\Pi }_{88}}  \\
			\end{matrix} \right]<0
\label{eq14f}
		\end{align}
		where
		\begin{align*}
			& {{\Pi }_{11}}={{P}_{a}}\tilde{A}+{{{\tilde{A}}}^{T}}{{P}_{a}}+\frac{{{\left\| {\tilde{A}} \right\|}^{2}}}{4{{\varepsilon }_{1}}}+\bar{d}^{2}{{(C\tilde{A})}^{T}}{{Q}_{3}}C\tilde{A} \cr
			& {{\Pi }_{14}}=-\bar{d}^{2}{{(C\tilde{A})}^{T}}{{Q}_{3}}B\tilde{L}(t) \cr
			& {{\Pi }_{17}}=-\bar{d}^{2}{{(C\tilde{A})}^{T}}{{Q}_{3}}(CX)\alpha \beta \tilde{L}(t) \cr
			& {{\Pi }_{18}}=-\bar{d}^{2}{{(C\tilde{A})}^{T}}{{Q}_{3}}, \,\, {{\Pi }_{22}}={{P}_{1}}+{{P}_{2}}-{{P}_{3}} \cr
			& {{\Pi }_{24}}=-{{P}_{1}}B\tilde{L}(t)+{{P}_{3}} ,\,\, {{\Pi }_{28}}=-{{P}_{1}}
		\end{align*}
		\begin{align*}
			 {{\Pi }_{33}}&=-{{P}_{1}}-{{P}_{3}} ,\,\,\, {{\Pi }_{34}}={{P}_{3}} \cr
			 {{\Pi }_{44}}&=-(1-{\varpi }){{P}_{3}}+\bar{d}^{2}{{(\beta \tilde{L}(t))}^{T}}{{P}_{4}}\beta \tilde{L}(t)\cr &-2{{P}_{4}}+\bar{d}^{2}{{(B\tilde{L}(t))}^{T}}{{Q}_{3}}B\tilde{L}(t) \cr
			 {{\Pi }_{45}}&=\frac{1}{2}(\pi \otimes {{I}_{q}})(-B\tilde{L}(t)) \cr
			 {{\Pi }_{47}}&=\bar{d}^{2}{{(-B\tilde{L}(t))}^{T}}{{Q}_{3}}CX(\alpha \beta \tilde{L}(t)) \cr
			 {{\Pi }_{48}}&=\bar{d}^{2}{{(B\tilde{L}(t))}^{T}}{{P}_{4}}+\bar{d}^{2}{{(B\tilde{L}(t))}^{T}}{{Q}_{3}}
		\end{align*}
		\begin{align*}
			& {{\Pi }_{55}}={{Q}_{1}}+{{Q}_{2}}+{{l}_{\max }}I-{{Q}_{3}}+{{\varepsilon }_{1}}{{\left\| (\pi \otimes {{I}_{q}})C \right\|}^{2}} \cr
			& {{\Pi }_{57}}={{Q}_{3}}+(\pi \otimes {{I}_{q}})(-\alpha \beta \tilde{L}(t)) \cr
			& {{\Pi }_{58}}=-(\pi \otimes {{I}_{q}}) ,\,\, {{\Pi }_{66}}=-{{Q}_{3}}-{{Q}_{1}} , \,
\, {{\Pi }_{67}}={{Q}_{3}} \cr
			& {{\Pi }_{77}}=-(1-{\varpi }){{Q}_{2}}-2{{Q}_{3}}+\bar{d}^{2}{{(\alpha \beta \tilde{L}(t))}^{T}}{{Q}_{3}}{{(\alpha \beta \tilde{L}(t))}} \cr
			& {{\Pi }_{78}}=(1/2)\bar{d}^{2}{{(\alpha \beta \tilde{L}(t))}^{T}}{{Q}_{3}} \cr
			& {{\Pi }_{88}}=-I+\bar{d}^{2}{{P}_{4}}+\bar{d}^{2}{{Q}_{3}}
		\end{align*}
\label{thm1}
\end{theorem}
\textbf{Proof}. A major factor contributing to stability Theorem \ref{thm1} is the convergence of $({{x}}(t),{{\eta }}(t),{{z }}(t))$ to the equilibrium point.  As we shall see, the output ${{y}_{i}}(t)$ at the equilibrium point is an optimal solution of  (\ref{eq4}). \\
We recall the system (\ref{eq13}), and we can obtain the equilibrium point $(\bar{x},\bar{\eta },\bar{z})$ from
\begin{align}
	& \dot{\bar{x}}(t)=0\Rightarrow (A-BK)(\tilde{x}(t)-X\tilde{\eta }(t))+X\tilde{\eta }(t)=0, \label{eq18}\\
	& \dot{\bar{\eta }}(t)=0\Rightarrow -\nabla \tilde{f}(y(t))-\alpha \beta \tilde{L}(t)\tilde{y}(t-d(t))-\beta \tilde{L}(t)\cr & \,\,\,\,\,\,\,\,\,\,\,\,\,\,\,\,\,\,\,\,\,\,\,\,\,\,\,\ \times(\tilde{\eta }(t-d(t))+\tilde{v}(t-d(t)))=0, \label{eq19}\\
	& \dot{\bar{z}}(t)=0\Rightarrow \alpha \beta \tilde{L}(t)\tilde{y}(t-d(t))=\alpha \beta \tilde{L}(t)C\bar{x}(t-d(t))=0. \label{eq20}
\end{align}
We assumed that $\Omega ([-\bar{d},0],{{\mathbb{R}}^{N{{n}_{i}}}})$ be the Banach space of continuous functions mapping the interval $[-\bar{d},0]$ into ${{\mathbb{R}}^{N}}$ with norm $\left\| \Phi  \right\|={{\sup }_{-\bar{d}\le y\le 0}}\left\| \Phi (y) \right\|$. The initial value can be considered by  $(y(0),\Phi (0))$ for delayed system (\ref{eq18})-(\ref{eq20}) in which $\Phi (0)$ is an element of $\Omega ([-\bar{d},0],{{\mathbb{R}}^{N{{n}_{i}}}})$. Since time delay appears in the $(y(t), \eta (t),z(t))$  component, thus, the initial value for $(y(t),\eta (t),z(t))$ is a function in $\Omega ([-\bar{d},0],{{\mathbb{R}}^{N{{n}_{i}}}})$.
\\
The equilibrium point is the solution, as we demonstrate in the sequel. Taking (\ref{eq18})-(\ref{eq20}) as a starting point, the equilibrium point meets
\begin{align}
	& (A-BK)(\bar{x}-X\bar{\eta })=0, \\
	& \sum\limits_{i=1}^{N}{\nabla {{{\tilde{f}}}_{i}}({{{\bar{y}}}_{i}})}=0, \\
	& {{{\bar{y}}}_{i}}={{C}_{i}}{{{\bar{x}}}_{i}}={{y}^{*}},
\end{align}
where ${{\text{y}}^{*}}\in {{\mathbb{R}}^{q}}$, and since $\tilde{f}$ is strongly convex, $\sum\limits_{i=1}^{N}{\nabla {{{\tilde{f}}}_{i}}({{y}^{*}})}=0,$ implies that ${{y}^{*}}$ is the optimal solution. Since $K$ is selected such that $(A-BK)$ is Hurwitz, it follows from (\ref{eq12}) and (\ref{eq13}) that
\begin{align*}
	& \bar{x}=X\bar{\eta }, \cr
	& \bar{y}=\bar{\eta }=1\otimes {{y}^{*}}
\end{align*}
Thus, $(\bar{x},\bar{\eta },\bar{z})=(X\bar{\eta },\bar{\eta },\bar{z})$ , and ${{y}_{i}}$ at the equilibrium point is the optimal solution of (\ref{eq4}).\\
Let $F(t)=\nabla \tilde{f}(\tilde{y}(t)+\bar{y}(t))-\nabla \tilde{f}(\bar{y}(t))$ and $\tilde{y}(t)=C\tilde{x}(t)$. Since $\tilde{x}(t)=x(t)-\bar{x}(t),$ $\tilde{\eta }(t)=\eta (t) -\bar{\eta }(t),$$\tilde{z}(t)=z(t)-\bar{z}(t),$ we further have the following system
\begin{align}
	 \dot{\bar{x}}(t)=&(A-BK)(\tilde{x}(t)-X\tilde{\eta }(t))+X\tilde{\eta }(t), \label{eq27} \\
	 \dot{\bar{\eta }}(t)=&F(t)-\alpha \beta \tilde{L}(t)\tilde{y}(t-d(t))-\beta \tilde{L}(t)\Big(\tilde{\eta }(t-d(t)) \cr &+\tilde{v}(t-d(t))\Big), \label{eq28} \\
	 \dot{\bar{z}}(t)=&\alpha \beta \tilde{L}(t)\tilde{y}(t-d(t)).  \label{eq299}
\end{align}
We consider the following candidate as a nonnegative Lyapunov function:
\begin{align}
	{{V}_{1}}(t)={{(\tilde{x}(t)-X\tilde{\eta }(t))}^{T}}{{P}_{a}}(\tilde{x}(t)-X\tilde{\eta }(t)),
	\label{eq29}
\end{align}
\begin{align}
	 {{V}_{2}}(t)=&0.5{{(\tilde{\eta }(t)+\tilde{v}(t))}^{T}}{{P}_{1}}(\tilde{\eta }(t)+\tilde{v}(t)) \cr
	& +\int_{t-\bar{d}}^{t}{{{(\tilde{\eta }(s)+\tilde{v}(s))}^{T}}{{P}_{2}}(\tilde{\eta }(s)+\tilde{v}(s))ds}\cr &+\int_{t-d(t)}^{t}{{{(\tilde{\eta }(s)+\tilde{v}(s))}^{T}}{{P}_{3}}(\tilde{\eta }(s)+\tilde{v}(s))ds} \cr
	&+\,\bar{d}\int_{-\bar{d}}^{0}{\int_{t+\theta }^{t}{{{(\tilde{\eta }(\theta )+\tilde{v}(\theta ))}^{T}}{{P}_{4}}(\tilde{\eta }(\theta )+\tilde{v}(\theta ))d\theta ds}},  \label{eq31}
\end{align}
\begin{align}
	 {{V}_{3}}(t)=&0.5{{{\tilde{y}}}^{T}}(t)(\Pi \otimes {{I}_{q}})\tilde{y}(t) \cr
	&+\int_{t-\bar{d}}^{t}{{{{\tilde{y}}}^{T}}(s){{Q}_{1}}\tilde{y}(s)ds}+\int_{t-d(t)}^{t}{{{{\tilde{y}}}^{T}}(s){{Q}_{2}}\tilde{y}(s)ds} \cr
	& +\,\bar{d}\int_{-\bar{d}}^{0}{\int_{t+\theta }^{t}{{{{\tilde{y}}}^{T}}(\theta ){{Q}_{3}}\tilde{y}(\theta )d\theta ds}}, \label{eq32}
\end{align}
where ${{P}_{a}},{{P}_{1}},{{P}_{2}},{{P}_{3}},{{P}_{4}},{{Q}_{1}},{{Q}_{2}},{{Q}_{3}}>0$ and $\Pi =diag\left\{ \pi  \right\}$ is defined in Lemma \ref{lem1}.
Therefore, the time derivative of (\ref{eq29}) along (\ref{eq27}) is characterized by
\begin{align}
	{{\dot{V}}_{1}}(t)={{(\tilde{x}(t)-X\tilde{\eta }(t))}^{T}}({{P}_{a}}{{\tilde{A}}^{T}}+\tilde{A}{{P}_{a}})(\tilde{x}(t)-X\tilde{\eta }(t)),
\label{eq27f}
\end{align}
where $\tilde{A}=A-BK$ is Hurwitz so that there is a positive define matrix ${\mathrm O}$ satisfying ${{P}_{a}}{{\tilde{A}}^{T}}+\tilde{A}{{P}_{a}}=-{\mathrm O}$ and ${{\varepsilon }_{1}}={{\lambda }_{mn}}({\mathrm O})>0$ .
Then, it follows from (\ref{eq31}) and (\ref{eq27})-(\ref{eq299}) that
\begin{align}
	 {{{\dot{V}}}_{2}}(t)=&{{\Big(\tilde{\eta }(t)+\tilde{v}(t)\Big)}^{T}}{{P}_{1}}\Big(-F(t)-\beta \tilde{L}(t)\big(\tilde{\eta }(t-d(t)) \cr
	& +\tilde{v}(t-d(t))\big)\Big)+{{\Big(\tilde{\eta }(t)+\tilde{v}(t)\Big)}^{T}}{{P}_{2}}\Big(\tilde{\eta }(t)+\tilde{v}(t)\Big) \cr
	& -{{\Big(\tilde{\eta }(t-\bar{d})+\tilde{v}(t-\bar{d})\Big)}^{T}}{{P}_{2}}\Big(\tilde{\eta }(t-\bar{d})+\tilde{v}(t-\bar{d})\Big) \cr
	& +{{\Big(\tilde{\eta }(t)+\tilde{v}(t)\Big)}^{T}}{{P}_{3}}\Big(\tilde{\eta }(t)+\tilde{v}(t)\Big)-\Big(1-\dot{d}(t)\Big)\Big(\tilde{\eta }(t-d(t)) \cr
	& +\tilde{v}(t-d(t)){{\Big)}^{T}}{{P}_{3}}\Big(\tilde{\eta }(t-d(t))+\tilde{v}(t-d(t))\Big) \cr
	& +{{{\bar{d}}}^{2}}{{\Big(\dot{\tilde{\eta }}(t)+\dot{\tilde{v}}(t)\Big)}^{T}}{{P}_{4}}\Big(\dot{\tilde{\eta }}(t)+\dot{\tilde{v}}(t)\Big) \cr
	& -\bar{d}\int_{t-d(t)}^{t}{{{(\dot{\tilde{\eta }}(s)+\dot{\tilde{v}}(s))}^{T}}{{P}_{4}}(\dot{\tilde{\eta }}(s)+\dot{\tilde{v}}(s))ds}. \label{eq33}
\end{align}
In (\ref{eq33}), ${{\bar{d}}^{2}}{{(\dot{\tilde{\eta }}(t)+\dot{\tilde{v}}(t))}^{T}}{{P}_{4}}(\dot{\tilde{\eta }}(t)+\dot{\tilde{v}}(t))$  can be replaced by the following relationships:
\begin{align}
	& {{{\bar{d}}}^{2}}\Big(-F(t)-\alpha \beta \tilde{L}(t)\tilde{y}(t-d(t))-\beta \tilde{L}(t)(\tilde{\eta }(t-d(t))\cr&+\tilde{v}(t-d(t)))+\alpha \beta \tilde{L}(t)\tilde{y}(t-d(t))\Big)^{T}{{P}_{4}}
	 \Big(-F(t)\cr &-\alpha \beta \tilde{L}(t)\tilde{y}(t-d(t))-\beta \tilde{L}(t)(\tilde{\eta }(t-d(t))+\tilde{v}(t-d(t)))\cr &+\alpha \beta \tilde{L}(t)\tilde{y}(t-d(t))\Big) .
\end{align}
According to Lemma \ref{lem11}, for last term in (\ref{eq33}) also we can have:
\begin{align}
	 -\bar{d}\int_{t-d(t)}^{t}&{{{(\dot{\tilde{\eta }}(s)+\dot{\tilde{v}}(s))}^{T}}{{P}_{4}}(\dot{\tilde{\eta }}(s)+\dot{\tilde{v}}(s))ds}\le  \cr
	& {{\left[ \begin{matrix}
				\tilde{\eta }(t)+\tilde{v}(t)  \\
				\tilde{\eta }(t-\bar{d})+\tilde{v}(t-\bar{d})  \\
				\tilde{\eta }(t-d(t))+\tilde{v}(t-d(t))  \\
			\end{matrix} \right]}^{T}}\left[ \begin{matrix}
		-{{P}_{4}} & {0} & {{P}_{4}}  \\
		{0} & -{{P}_{4}} & {{P}_{4}}  \\
		{0} & {0} & -2{{P}_{4}}  \\
	\end{matrix} \right]  \cr &  \times \left[ \begin{matrix}
		\tilde{\eta }(t)+\tilde{v}(t)  \\
		\tilde{\eta }(t-\bar{d})+\tilde{v}(t-\bar{d})  \\
		\tilde{\eta }(t-d(t))+\tilde{v}(t-d(t))  \\
	\end{matrix} \right].  \label{eq35}
\end{align}
Meanwhile, the time derivative of the first term in (\ref{eq32}) is expressed as
\begin{align}
	 {{{\tilde{y}}}^{T}}(t)(\Pi \otimes & {{I}_{q}})\dot{\tilde{y}}(t)={{{\tilde{y}}}^{T}}(t)(\Pi \otimes {{I}_{q}})C\dot{\tilde{x}}(t) \cr
	=&{{{\tilde{y}}}^{T}}(t)(\Pi \otimes {{I}_{q}})C\Big((A-BK)(\tilde{x}(t)-X\tilde{\eta }(t))+X\dot{\tilde{\eta }}(t)\Big) \cr
	 =&{{{\tilde{y}}}^{T}}(t)(\Pi \otimes {{I}_{q}})C\Big((A-BK)(\tilde{x}(t)-X\tilde{\eta }(t))\Big)\cr &+{{{\tilde{y}}}^{T}}(t)(\Pi \otimes {{I}_{q}})CX\Big(-F(t)-\alpha \beta \tilde{L}(t)\tilde{y}(t-d(t))\cr & - \beta \tilde{L}(t)\big(\tilde{\eta }(t-d(t))+\tilde{v}(t-d(t))\big)\Big),  \label{eq36}
\end{align}
for the first term in (\ref{eq36}), by using the stated concept in Lemma \ref{lem1} and the Young's inequality, i.e., ${{Y}^{T}}F\le {{Y}^{T}}Y+{{F}^{T}}F$ , one has:
\begin{align}
	& {{{\tilde{y}}}^{T}}(t)(\Pi \otimes {{I}_{q}})C\tilde{A}(\tilde{x}(t)-X\tilde{\eta }(t))\le {{\varepsilon }_{1}}{{\left\| (\Pi \otimes {{I}_{q}})C \right\|}^{2}}{{{\tilde{y}}}^{T}}(t)\tilde{y}(t) \cr
	& +\frac{{{\left\| {\tilde{A}} \right\|}^{2}}}{4{{\varepsilon }_{1}}}{{(\tilde{x}(t)-X\tilde{\eta }(t))}^{T}}(\tilde{x}(t)-X\tilde{\eta }(t)),
\end{align}
and for the fourth term in (\ref{eq36}) by using the Young inequality, we can have
\begin{align}
	& {{{\tilde{y}}}^{T}}(t)(\Pi \otimes {{I}_{q}})(-\beta \tilde{L}(t))\Big(\tilde{\eta }(t-d(t))+\tilde{v}(t-d(t))\Big)= \cr
	& \,\,{{{\tilde{y}}}^{T}}(t)\frac{1}{2}(\Pi \otimes {{I}_{q}})(-\beta \tilde{L}(t))\Big(\tilde{\eta }(t-d(t))+\tilde{v}(t-d(t))\Big) \cr
	& \,\,+ {{\Big(\tilde{\eta }(t-d(t))+\tilde{v}(t-d(t))\Big)}^{T}}\frac{1}{2}(\Pi \otimes {{I}_{q}})(-\beta \tilde{L}(t))\tilde{y}(t).
\end{align}
Also, for the second term in (\ref{eq32}), we can have
\begin{align}
	& \int_{t-\bar{d}}^{t}{{{{\tilde{y}}}^{T}}(s){{Q}_{1}}\tilde{y}(s)ds}={{{\tilde{y}}}^{T}}(t){{Q}_{1}}\tilde{y}(t)-{{{\tilde{y}}}^{T}}(t-\bar{d}){{Q}_{1}}\tilde{y}(t-\bar{d})
, \end{align}
and for third term in (\ref{eq32}), one has
\begin{align}
&\int_{t-d(t)}^{t}{{{{\tilde{y}}}^{T}}(s){{Q}_{2}}\tilde{y}(s)ds}	=\cr &{{{\tilde{y}}}^{T}}(t){{Q}_{2}}\tilde{y}(t)-(1-\dot{d}(t)){{{\tilde{y}}}^{T}}(t-d(t)){{Q}_{2}}\tilde{y}(t-d(t)),
\label{eq35f}
\end{align}
and for the fourth term in (\ref{eq32}), we can have
\begin{align}
&	\bar{d}\int_{-\bar{d}}^{0}{\int_{t+\theta }^{t}{{{{\tilde{y}}}^{T}}(\theta ){{Q}_{3}}\tilde{y}(\theta )d\theta ds}}=\cr & {{\bar{d}}^{2}}{{(\dot{\tilde{y}}(t))}^{T}}{{Q}_{3}}\dot{\tilde{y}}(t)-\bar{d}\int_{t-d(t)}^{t}{{{(\dot{\tilde{y}}(s))}^{T}}{{Q}_{3}}\dot{\tilde{y}}(s)ds}
\label{eq41}.
	\end{align}
The first term in equation (\ref{eq41}) can be convert to the following form:
\begin{align}
	 {{{\bar{d}}}^{2}}{{(\dot{\tilde{y}}(t))}^{T}}&{{Q}_{3}}\dot{\tilde{y}}(t)={{{\bar{d}}}^{2}}\Big({{(C\dot{\tilde{x}}(t))}^{T}}{{Q}_{3}}C\dot{\tilde{x}}(t)\Big) \cr
	& ={{{\bar{d}}}^{2}}{{\Big(C\big((A-BK)(\tilde{x}(t)-X\tilde{\eta }(t))+X\tilde{\eta }(t)\big)\Big)}^{T}}{{Q}_{3}} \cr
	&  \times \Big(C\big((A-BK)(\tilde{x}(t)-X\tilde{\eta }(t))+X\tilde{\eta }(t)\big)\Big) \cr
	& ={{{\bar{d}}}^{2}}\Big(C\Big((A-BK)(\tilde{x}(t)-X\tilde{\eta }(t)) \cr
	& +X\big(-F-\alpha \beta \tilde{L}(t)\tilde{y}(t-d(t))-\beta \tilde{L}(t)(\tilde{\eta }(t-d(t)) \cr
	& +\tilde{v}(t-d(t)))\big)\Big){{\Big)}^{T}}{{Q}_{3}}\Big(C\big((A-BK)(\tilde{x}(t)-X\tilde{\eta }(t)) \cr
	& +X(-F-\alpha \beta \tilde{L}(t)\tilde{y}(t-d(t))-\beta \tilde{L}(t)(\tilde{\eta }(t-d(t)) \cr
	& +\tilde{v}(t-d(t))))\big)\Big)={{{\bar{d}}}^{2}}\Big(C(A-BK)(\tilde{x}(t)-X\tilde{\eta }(t)) \cr
	& +(-F-\alpha \beta \tilde{L}(t)\tilde{y}(t-d(t))-\beta \tilde{L}(t)(\tilde{\eta }(t-d(t)) \cr
	& +\tilde{v}(t-d(t)))){{\Big)}^{T}}{{Q}_{3}}\Big(C(A-BK)(\tilde{x}(t)-X\tilde{\eta }(t)) \cr
	& +(-F-\alpha \beta \tilde{L}(t)\tilde{y}(t-d(t))-\beta \tilde{L}(t)(\tilde{\eta }(t-d(t)) \cr
	& +\tilde{v}(t-d(t))))\Big).
\label{eq37f}
\end{align}
According to Lemma \ref{lem11}, the second term in (\ref{eq41}) can be structured in the following way
\begin{align}
&-\bar{d}\int_{t-d(t)}^{t}{{{(\dot{\tilde{y}}(s))}^{T}}{{Q}_{3}}\dot{\tilde{y}}(s)ds} \le  \cr
	& {{\left[ \begin{matrix}
				\tilde{y }(t)  \\
				\tilde{y }(t-\bar{d})  \\
				\tilde{y }(t-d(t))  \\
			\end{matrix} \right]}^{T}}\left[ \begin{matrix}
		-{{Q}_{3}} & {0} & {{Q}_{3}}  \\
		{0} & -{{Q}_{3}} & {{Q}_{3}}  \\
		{0} & {0} & -2{{Q}_{3}}  \\
	\end{matrix} \right]   \left[ \begin{matrix}
		\tilde{y }(t)  \\
		\tilde{y }(t-\bar{d})  \\
		\tilde{y }(t-d(t))  \\
	\end{matrix} \right].  \label{eq38f}
\end{align}
We define the following variable
\begin{align}
	 \bar{\Phi} (t)=&\Big[{{(\tilde{x}-X\tilde{\eta })}^{T}},{{(\tilde{\eta }(t)+\tilde{v}(t))}^{T}},{{(\tilde{\eta }(t-\bar{d})+\tilde{v}(t-\bar{d}))}^{T}}, \cr
	& {{(\tilde{\eta }(t-d(t))+\tilde{v}(t-d(t)))}^{T}},{{{\tilde{y}}}^{T}}(t),{{{\tilde{y}}}^{T}}(t-\bar{d}),\cr&{{{\tilde{y}}}^{T}}(t-d(t)),{{F}^{T}(t)}\Big]^{T}
. \end{align}
According to the definition of Lipshitz condition, i.e., Assumption \ref{ass2}, we can have ${{l}_{\max }}{{\tilde{y}}^{T}}(t)\tilde{y}(t)-{{F}^{T}}(t)F(t)\ge 0$. This inequality then is added to $\dot{V}(t)={{\dot{V}}_{1}}(t)+{{\dot{V}}_{2}}(t)+{{\dot{V}}_{3}}(t)-{{l}_{\max }}{{\tilde{y}}^{T}}(t)\tilde{y}(t)-{{F}^{T}}(t)F(t)$,
and by considering  $(\ref{eq27f})-(\ref{eq35f})$, $(\ref{eq37f})-(\ref{eq38f})$, ${{C}_{i}}{{X}_{i}}={{I}_{q}}$ from (\ref{eq12}), $(1-\dot{d}(t))\le (1-\varpi )$,  the Schur complement and the variable $\Omega (t)$, LMI  (\ref{eq14f}) is obtained. If LMI (\ref{eq14f})  is satisfied, $\dot{V}(t)\le {\bar{\Phi }^{T}}(t)\Pi \bar\Phi (t) \le  0$ is kept. It completes the proof.\\
\begin{remark}
  By applying the Lyapunov-Krasovskii functional to time-delay systems \cite{fnf}, delay-dependent LMI conditions are derived. This condition provides a way to estimate the upper bound on the time-delay in an optimization problem (\ref{eq6})-(\ref{eq11}).
\end{remark}
\begin{remark}
  Notice that the feasibility of LMI (\ref{eq14f}) depends on the number of agents in a networked system. There would be difficulties in verifying the LMI conditions when there are a large number of agents. Consequently, this proposed approach has a major drawback regarding larger networked systems.
\end{remark}

\subsection{Fast Time-Varying and Without Delays Cases}
We designed a delay-dependent condition for the considered distributed optimization problem for heterogeneous linear multi-agent systems in the above.  It is therefore expected that time-varying delays would be slow in such circumstances. In other words, the time delay $\varpi$ must be less than one; for example, $\varpi \in \left[ 0,1 \right)$.  In the absence of such a constraint on time delays, the following delay-dependent conditions can be obtained without $\varpi$, in which time delays can vary at any speed. Bounded communications with unknown dynamics can also be handled by these conditions.
\begin{theorem}
  We assumed  Assumptions \ref{ass1}-\ref{ass2} are satisfied, under the introduced distributed optimization approach in (\ref{eq6})-(\ref{eq11}),  system (\ref{eq4}) with initial value in ${{x}_{i}}(0),{{\eta }_{i}}(0)$ and ${{z}_{i}}(0)$  convergent to the optimal solution for any delay $d(t)\in \left[ 0,\bar{d} \right]$, if ${{K}_{i}}$ is selected to that ${{A}_{i}}-{{B}_{i}}{{K}_{i}}$ is Hurwitz, ${{U}_{i}},\,{{W}_{i}},\,{{X}_{i}}$ are solutions to (\ref{eq12}), if all the LMI conditions in Theorem \ref{thm1} with $ {{P}_{3}}=0$ and ${{Q}_{2}}=0$.
\end{theorem}
\textbf{Proof:}In the considered Lyapunov functions (\ref{eq29})-(\ref{eq32}), with $ \varpi $ only has relationships with the following terms according to Theorem \ref{thm1}.
\begin{align*}
 & \int_{t-d(t)}^{t}{{{{\tilde{y}}}^{T}}(s){{Q}_{2}}\tilde{y}(s)ds}  \cr &
\int_{t-d(t)}^{t}{{{(\tilde{\eta }(s)+\tilde{v}(s))}^{T}}{{P}_{3}}(\tilde{\eta }(s)+\tilde{v}(s))ds}
\end{align*}
The following Lyapunov functional is taken into consideration:
\begin{align}
	{{V}_{1}}(t)={{(\tilde{x}(t)-X\tilde{\eta }(t))}^{T}}{{P}_{a}}(\tilde{x}(t)-X\tilde{\eta }(t)),
	\label{eq290}
\end{align}
\begin{align}
	 {{V}_{2}}(t)=&0.5{{(\tilde{\eta }(t)+\tilde{v}(t))}^{T}}{{P}_{1}}(\tilde{\eta }(t)+\tilde{v}(t)) \cr
	& +\int_{t-\bar{d}}^{t}{{{(\tilde{\eta }(s)+\tilde{v}(s))}^{T}}{{P}_{2}}(\tilde{\eta }(s)+\tilde{v}(s))ds}\cr
	&+\,\bar{d}\int_{-\bar{d}}^{0}{\int_{t+\theta }^{t}{{{(\tilde{\eta }(\theta )+\tilde{v}(\theta ))}^{T}}{{P}_{4}}(\tilde{\eta }(\theta )+\tilde{v}(\theta ))d\theta ds}},  \label{eq310}
\end{align}
\begin{align}
	 {{V}_{3}}(t)=&0.5{{{\tilde{y}}}^{T}}(t)(\Pi \otimes {{I}_{q}})\tilde{y}(t) \cr
	&+\int_{t-\bar{d}}^{t}{{{{\tilde{y}}}^{T}}(s){{Q}_{1}}\tilde{y}(s)ds}\cr
	& +\,\bar{d}\int_{-\bar{d}}^{0}{\int_{t+\theta }^{t}{{{{\tilde{y}}}^{T}}(\theta ){{Q}_{3}}\tilde{y}(\theta )d\theta ds}}, \label{eq320}
\end{align}
The conclusion can be reached directly by following the same procedure as in Theorem \ref{thm1}. \quad \\

If there is no time delay, sufficient condition in form of LMI can be obtained to guarantee that the optimization problem of the delay-free system (\ref{eq6})-(\ref{eq9}) will be convergent. Lyapunov's function can be described as follows:
\begin{align}
	& {{V}_{1}}(t)={{(\tilde{x}(t)-X\tilde{\eta }(t))}^{T}}{{P}_{a}}(\tilde{x}(t)-X\tilde{\eta }(t)), \cr &
	 {{V}_{2}}(t)=0.5{{(\tilde{\eta }(t)+\tilde{v}(t))}^{T}}{{P}_{1}}(\tilde{\eta }(t)+\tilde{v}(t)) \cr &
	 {{V}_{3}}(t)=0.5{{{\tilde{y}}}^{T}}(t)(\Pi \otimes {{I}_{q}})\tilde{y}(t)
\end{align}
then, in this case, we can obtain delay-free sufficient conditions using the same procedures as in Theorem \ref{thm1}.
\section{Simulation Results}
\label{sec_sr}
\begin{figure}[h]
\centerline{\includegraphics[width=3.5in, height=3.0 in]{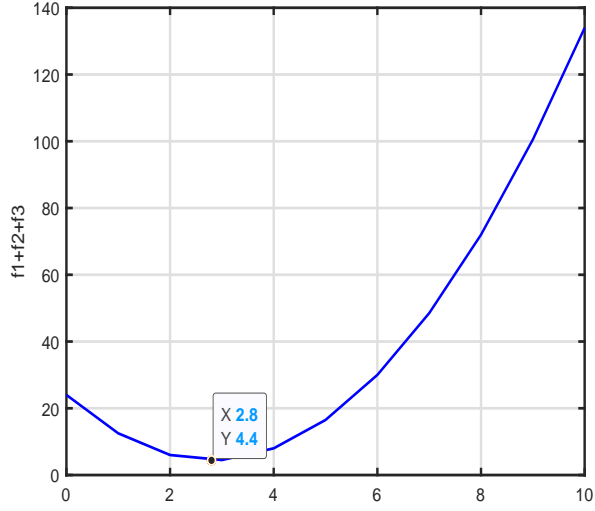}}
\caption{Images of the values of ${f}_{1}+{f}_{2}+{f}_{3}$.}
\label{fig_3}
\end{figure}

\begin{figure}[h]
\centerline{\includegraphics[width=3.0in, height=2.0 in]{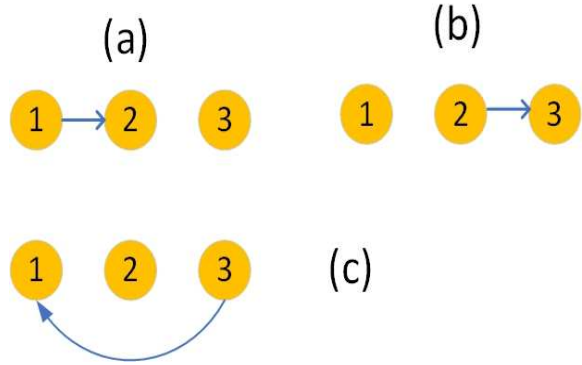}}
\caption{A network of three heterogeneous agents and their communication digraphs.}
\label{fig_1}
\end{figure}

\begin{figure}[h]
\centerline{\includegraphics[width=3.5in, height=3.0 in]{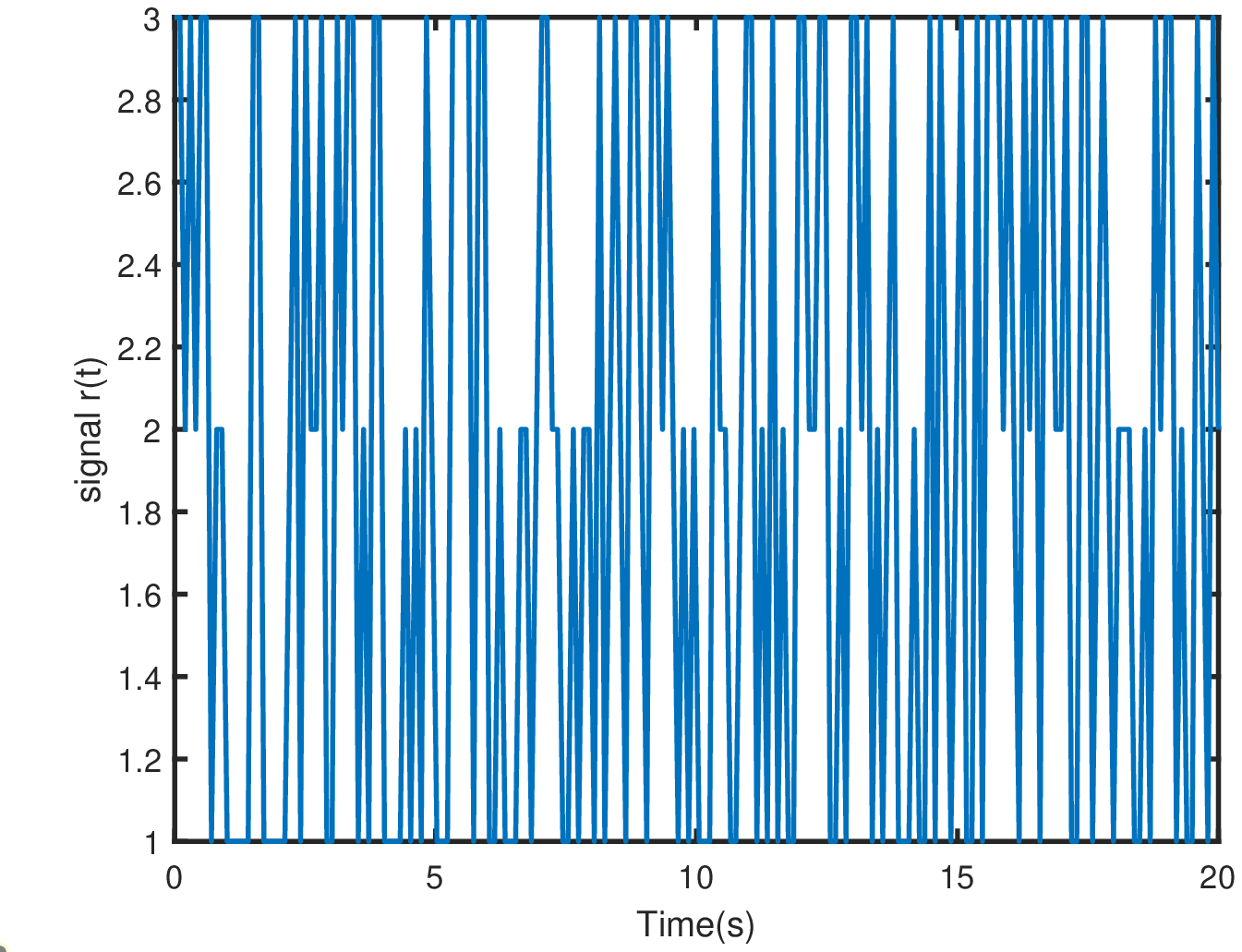}}
\caption{The signal $r(t)$ for the considered random digraphs.}
\label{fig_2}
\end{figure}

\begin{figure}[h]
\centerline{\includegraphics[width=3.5in, height=3.0 in]{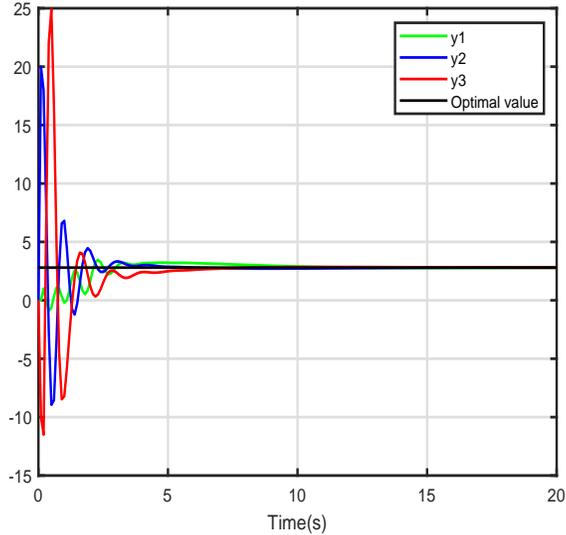}}
\caption{The transient behavior of the output ${y}_{i}$ for each agent with $\alpha =1,\,\,\beta =0.75,\,\,d(t)=0.1s$.}
\label{fig_33}
\end{figure}

\begin{figure}[h]
\centerline{\includegraphics[width=3.5in, height=3.0 in]{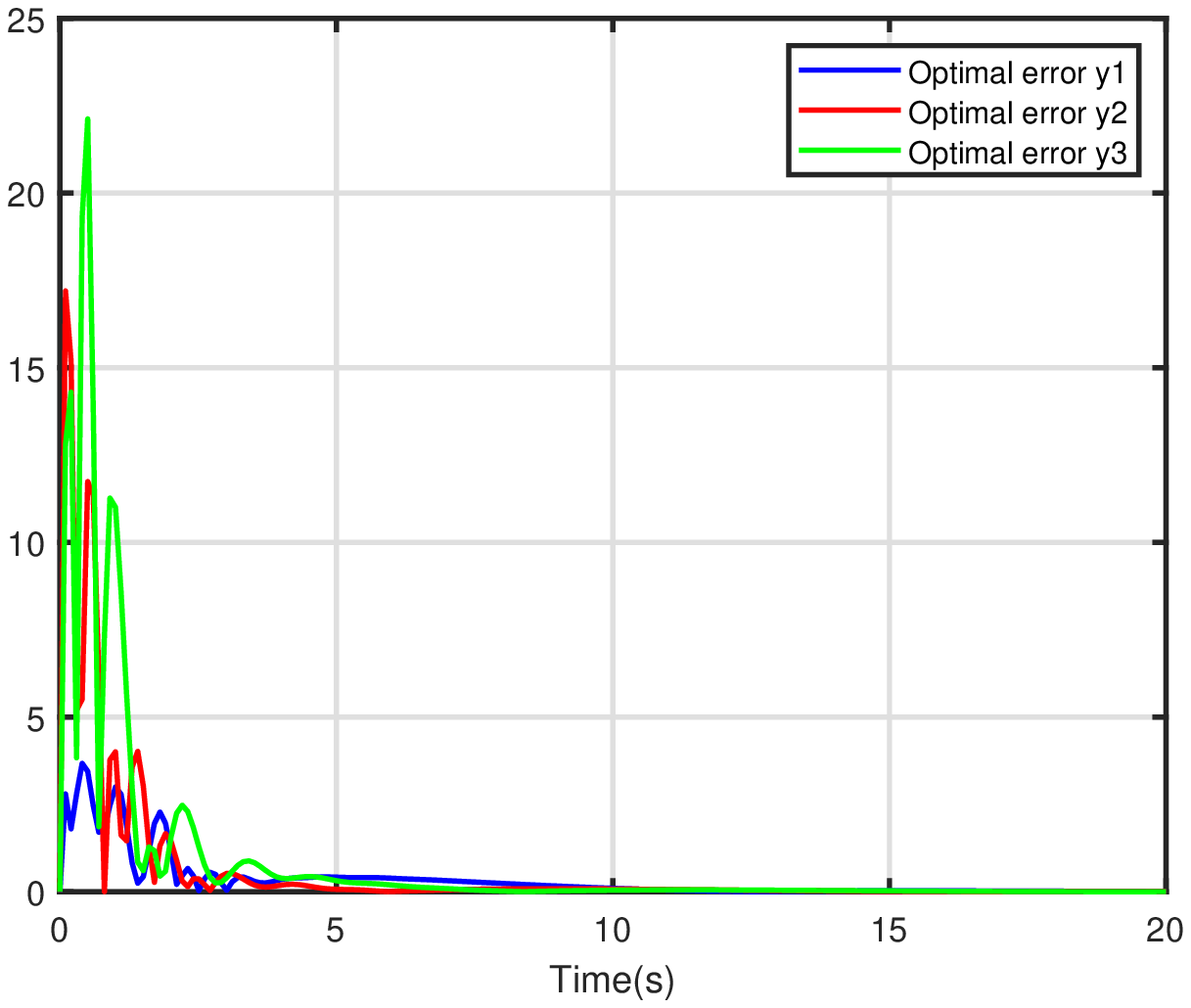}}
\caption{The optimal tracking error for each agent with $\alpha =1,\,\,\beta =0.75,\,\,d(t)=0.1s$.}
\label{fig_44}
\end{figure}

\begin{figure}[h!]
\centerline{\includegraphics[width=3.5in, height=3.0 in]{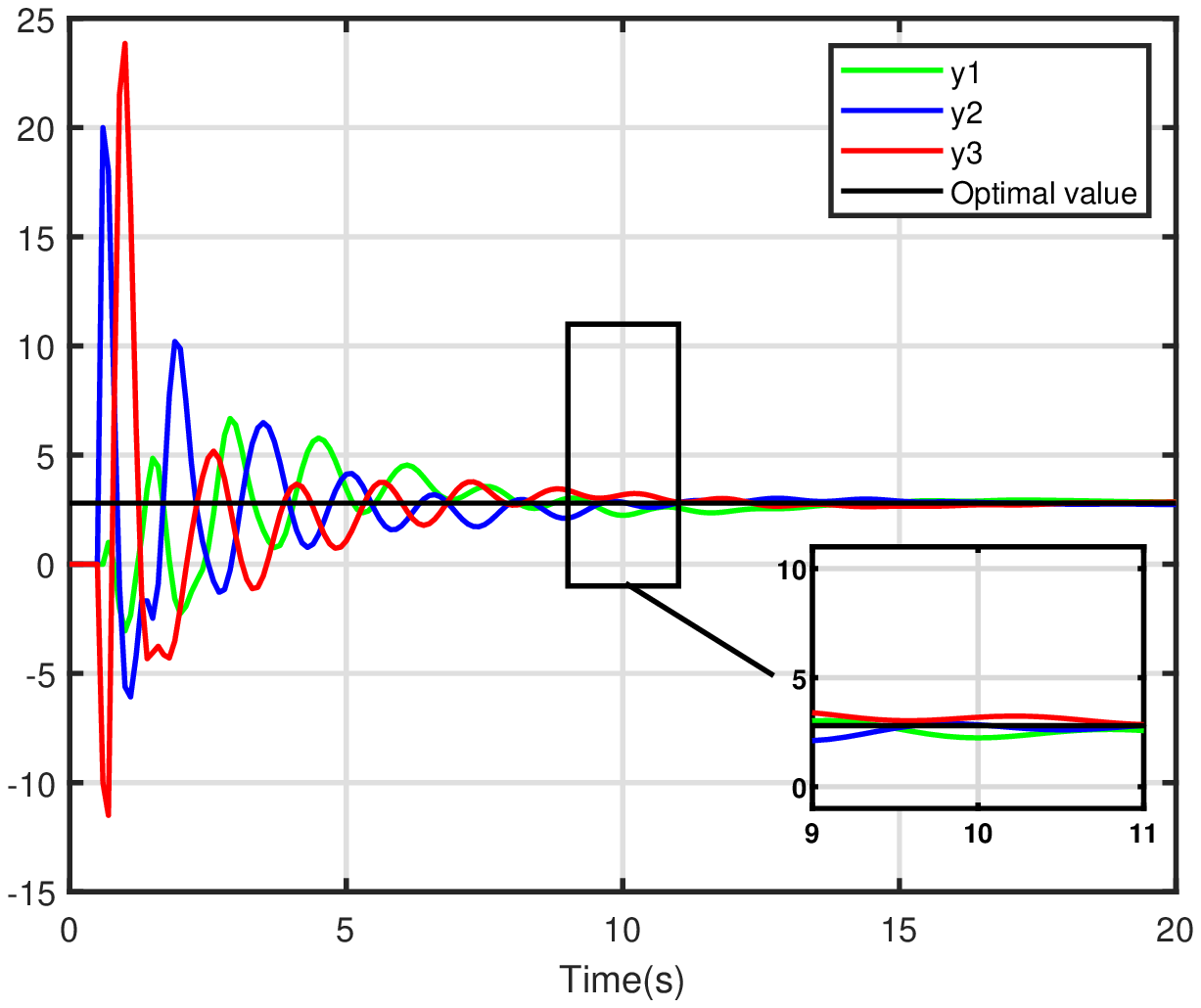}}
\caption{The transient behavior of the output ${y}_{i}$ for each agent with $\alpha =1,\,\,\beta =0.75,\,\,d(t)=0.4s$.}
\label{fig_6}
\end{figure}

\begin{figure}[h!]
\centerline{\includegraphics[width=3.5in, height=3.0 in]{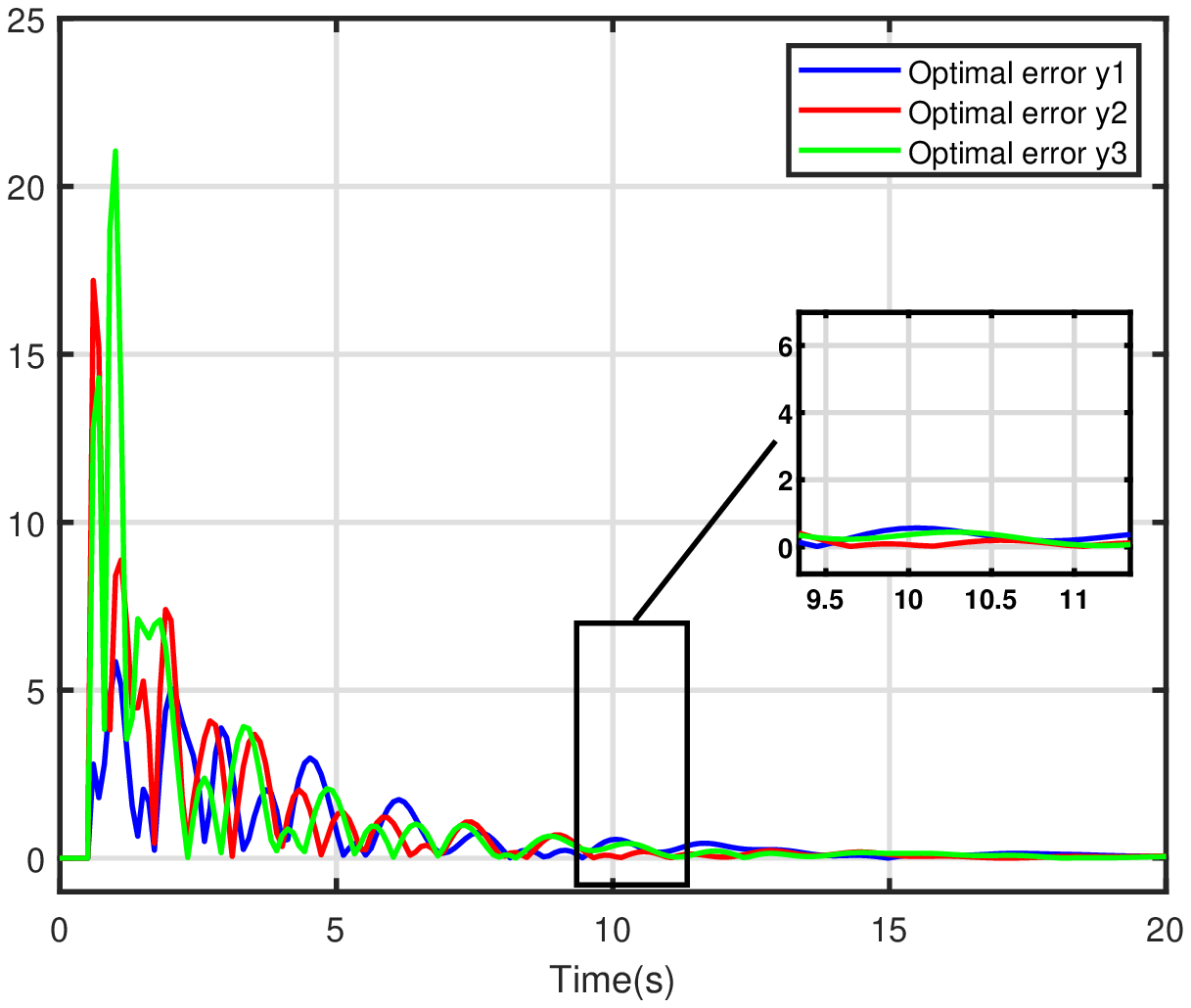}}
\caption{The optimal tracking error for each agent with $\alpha =1,\,\,\beta =0.75,\,\,d(t)=0.4s$.}
\label{fig_7}
\end{figure}

To validate the effectiveness of the proposed approach, numerical simulation examples are provided in this section. We considered a network of heterogeneous linear systems which includes three agents  with different dimensions and  the dynamics of them  is described as follows:
\begin{align*}
  & i=1\to {{A}_{1}}=\left[ \begin{matrix}
   0 & 1  \\
   0 & 0  \\
\end{matrix} \right],\,\,{{B}_{1}}=\left[ \begin{matrix}
   0 & 1  \\
   1 & -2  \\
\end{matrix} \right],\,\,{{C}_{1}}=\left[ \begin{matrix}
   1 & 1  \\
\end{matrix} \right] \cr
 & i=2\to {{A}_{2}}=\left[ \begin{matrix}
   0 & -1  \\
   1 & -2  \\
\end{matrix} \right],\,\,{{B}_{2}}=\left[ \begin{matrix}
   1 & 0  \\
   3 & -1  \\
\end{matrix} \right],\,\,{{C}_{2}}=\left[ \begin{matrix}
   -1 & 1  \\
\end{matrix} \right] \cr
 & i=3\to {{A}_{3}}=\left[ \begin{array}{*{35}{l}}
   0 & 1 & 0  \\
   0 & 0 & 1  \\
   0.5 & 1 & -2  \\
\end{array} \right],\,\,{{B}_{3}}=\left[ \begin{array}{*{35}{l}}
   1 & 0  \\
   0 & 1  \\
   1 & 0  \\
\end{array} \right], {{C}_{3}}=\left[ \begin{matrix}
   1 & -1 & 1  \\
\end{matrix} \right].
\end{align*}
Our first step should be to verify Assumption \ref{asp1}, for example, for agent $i=1$, we have: $rank(\left[ \begin{matrix}
   {{C}_{1}}{{B}_{1}} & {{0}_{1\times 2}}  \\
   -{{A}_{1}}{{B}_{1}} & {{B}_{1}}  \\
\end{matrix} \right])=3={{n}_{1}}+{{q}_{1}},\,\,\,{{n}_{1}}=2,\,\,{{q}_{1}}=1$.
We can use the above relationship to conclude that Assumption \ref{asp1} has been satisfied, and the next steps of the proposed approach may now be implemented.
Following is an optimization problem that we considered:
\begin{align}
F(\theta )=\sum\limits_{i=1}^{N=3}{{{f}_{i}}(\theta )},
\label{eq446}
\end{align}
in which the local objective functions of each agent are defined as follows:
\begin{align*}
  & {{f}_{1}}(\theta )=\frac{1}{2}{{\theta }^{2}}-1,\,\,{{f}_{2}}(\theta )={{(\theta -4)}^{2}}, \\
 & {{f}_{3}}(\theta )={{(\theta -3)}^{2}},
\end{align*}
it is apparent that all the local objective functions are strongly convex and Lipschitz global.
It can be simply verified that the general cost function (\ref{eq446}), gets its minimum value ${{F}_{\min }}=4.4$ at ${{x}_{\min }}=2.8$ as shown in Fig. \ref{fig_3}. In the first simulation, the delay value is supposed to be constant $d(t)=0.1s$.
In what follows, the distributed optimization approach (\ref{eq6})-(\ref{eq11}) is executed for a networked of three agents where they exchange their information with each other according to digraphs in Fig. \ref{fig_1}. In Fig. \ref{fig_2}, $r(t)$ represents the communication topology between agents, which can switch between three scenarios.
Fig. \ref{fig_1} shows three randomly generated digraphs for a team of three agents, in which the random process has been captured by the Markov chain, with $S=\{1,2,3\}$ as the state space;
it has the transition rate matrix as follows:
\begin{align*}
\Upsilon =\left( \begin{array}{*{35}{l}}
   -0.2 & 0.1 & 0.1  \\
   0.2 & -0.6 & 0.4  \\
   0.02 & 0.08 & -0.1  \\
\end{array} \right)
\end{align*}
in which all off-diagonal elements are nonnegative and the row summation is zero. $\{0.4772, 0.2612,0.3235\}$ are the initial values of the distribution. It is evident from the diagram that the graphs are disconnected, but that the union of the graphs contains a directed spanning tree to meet Assumption \ref{ass6}. For these three agents, by choosing ${{K}_{i}}$ as the controller gain matrix, ${{A}_{i}}-{{B}_{i}}{{K}_{i}}$ can be represented as a Hurwitz matrix, which is presented as follows:
\begin{align*}
  & {{K}_{1}}=\left[ \begin{array}{*{35}{l}}
   8 & 7  \\
   4 & 1  \\
\end{array} \right],\,\,{{K}_{2}}=\left[ \begin{array}{*{35}{l}}
   3 & -1  \\
   8 & -5  \\
\end{array} \right], \cr
 & {{K}_{3}}=\left[ \begin{array}{*{35}{l}}
   6.3333 & 1 & -1.3333  \\
   0 & 2 & 1  \\
\end{array} \right]
\end{align*}
and the required gain matrices based on (\ref{eq12}) are obtained as follows:
\begin{align*}
  & {{U}_{1}}={{\left[ \begin{matrix}
   1 & 0.5  \\
\end{matrix} \right]}^{T}},\,\,{{W}_{1}}={{\left[ \begin{matrix}
   1.5 & 0.5  \\
\end{matrix} \right]}^{T}},\,\,{{X}_{1}}={{\left[ \begin{matrix}
   0.5 & 0.5  \\
\end{matrix} \right]}^{T}}, \cr
 & {{U}_{2}}={{\left[ \begin{matrix}
   -0.5 & 0  \\
\end{matrix} \right]}^{T}},\,\,{{W}_{2}}={{\left[ \begin{matrix}
   -0.5 & -2  \\
\end{matrix} \right]}^{T}},\,\,{{X}_{2}}={{\left[ \begin{matrix}
   -0.5 & 0.5  \\
\end{matrix} \right]}^{T}}, \cr
 & {{U}_{3}}={{\left[ \begin{matrix}
   -1 & 0  \\
\end{matrix} \right]}^{T}},\,\,{{W}_{3}}={{\left[ \begin{matrix}
   0 & -1  \\
\end{matrix} \right]}^{T}},\,\,{{X}_{3}}={{\left[ \begin{matrix}
   0 & -1 & 0  \\
\end{matrix} \right]}^{T}}.
\end{align*}
We assume that $d(t)$ is a constant communication delay, i.e., $\bar{d}=0$. This set of parameters allows LMI (\ref{eq14f}) in Theorem \ref{thm1} to be solved. In Figs.  \ref{fig_33}-\ref{fig_44} we can see the simulation results. It can be seen that all the agent’s outputs converge to the optimal solution 2.8.  Fig. \ref{fig_44}  displays the error of tracking the optimal point which demonstrates that the convergence performance is satisfactory. To investigate the impact of increasing communication delays on the proposed distributed optimization approach, simulation results are provided in Figs. \ref{fig_6}-\ref{fig_7}.


\begin{figure}[h!]
\centerline{\includegraphics[width=3.5in, height=3.0 in]{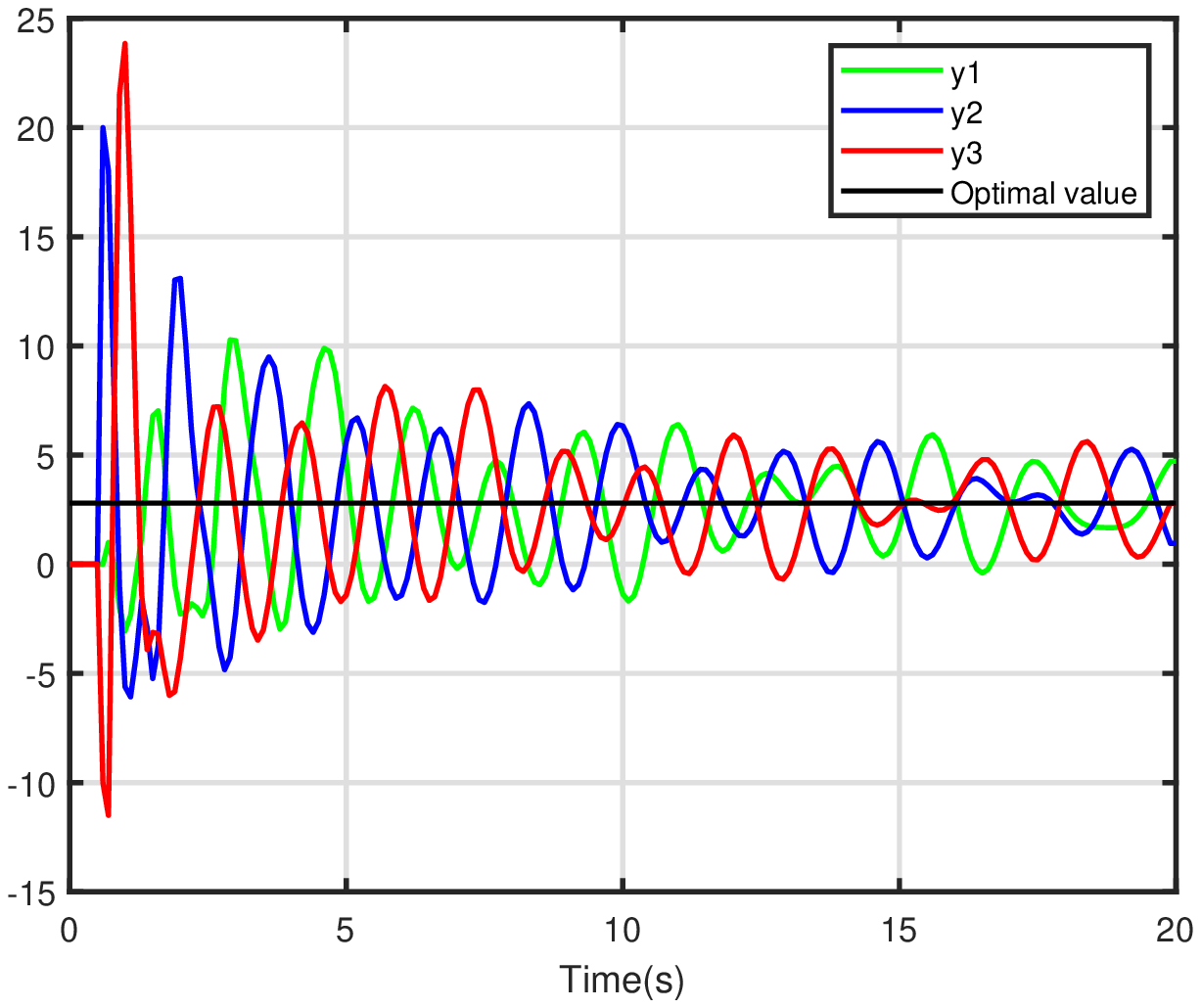}}
\caption{The transient behavior of the output ${y}_{i}$ for each agent with $\alpha =1,\,\,\beta =0.75,\,\,d(t)=0.7s$.}
\label{fig_8}
\end{figure}

\begin{figure}[h!]
\centerline{\includegraphics[width=3.5in, height=3.0 in]{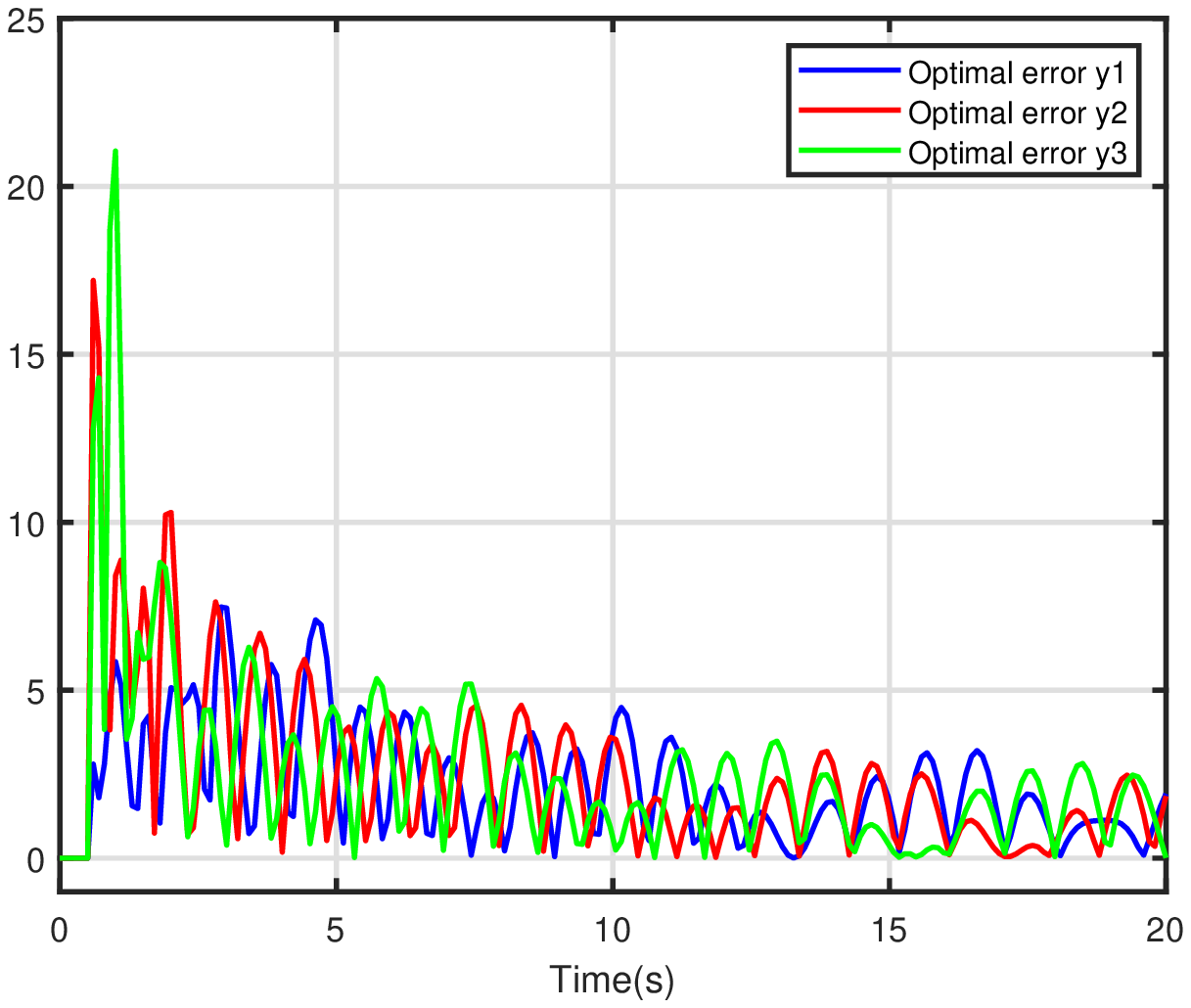}}
\caption{The optimal tracking error for each agent with $\alpha =1,\,\,\beta =0.75,\,\,d(t)=0.7s$.}
\label{fig_9}
\end{figure}

As can be seen in Figs. \ref{fig_6}-\ref{fig_7}, unlike the performance of our optimization problem for $d(t)=0.1s$, the agents' outputs had not converged to the optimal point at 10s.
The simulation results for $d(t)=0.7s$ is shown in Figs. \ref{fig_8}-\ref{fig_9}, which shows that the outputs are not converging.
The maximum allowable delays for the introduced method are reported in Table \ref{tab1} to make evolution convenient for finding a relationship between increasing delay and convergence. In the case of a delay with $0.7s$, the proposed method can not meet the converging criterion.
\begin{table}
\label{tab1}
\caption{The performance of the proposed optimization approach}
\label{table}
\setlength{\tabcolsep}{3pt}
\begin{tabular}{|p{85pt}|p{75pt}|p{75pt}|}
\hline
Communication delay&
The feasibility of LMI (\ref{eq14f})&
Convergence \\
\hline
$d(t)=0.1s $&
feasible&
Yes \\
$d(t)=0.4s$&
feasible&
Yes \\
$d(t)=0.7s$&
Not feasible&
No \\
\hline
\end{tabular}
\label{tab1}
\end{table}


\section{Conclusion}
\label{sec_c}
We proposed a distributed optimization method for multi-agent systems in this paper.  In the addressed approach, each agent was assumed to have its own dynamic that could differ from another agents in the network. As well, communication delays between agents in the network have been taken into account. Also, taking into account the unreliability of networks, we assume that random digraphs illustrate the communication digraphs of a network of heterogeneous linear systems. Lastly, delay-dependent sufficient conditions are derived by means of linear matrix inequality to prove convergence to optimal solutions. In order to increase the usefulness of the proposed method for the practical systems in future studies, this work can be extended to include the following issues:
(1) Considering the quantization effects in exchanging information within a network \cite{fn1};
(2) Considering non-convex optimization problem. \cite{fn66}.

\section*{Conflict of interest} The authors declare that they have no conflict of interest.

\bibliographystyle{apacite}
\bibliography{ref1}

\end{document}